\documentclass[journal ]{new-aiaa}
\usepackage{textcomp}

\usepackage{graphicx}
\usepackage{amsmath}
\usepackage[version=4]{mhchem}
\usepackage{siunitx}
\usepackage{longtable,tabularx}
\title{Optimal Trajectories for Multiple-UAS Simultaneous Target Acquisition with Obstacle Avoidance}

\author{Michael D. Zollars\footnote{Lt Col, USAF, Assistant Professor, Department of Aeronautics and Astronautics, AIAA Member.} }
\affil{Air Force Institute of Technology, Wright-Patterson Air Force Base, OH, 45433}

\author{David J. Grymin\footnote{Controls Science Engineer, Controls Science Center, Air Force Research Laboratory, WPAFB, OH 45433, AIAA Senior Member.} and Isaac E. Weintraub\footnote{Electronics Engineer, Controls Science Center, Air Force Research Laboratory, WPAFB, OH 45433, AIAA Senior Member.}}
\affil{Air Force Research Laboratory, Wright-Patterson Air Force Base, Ohio, 45433}

\begin{document}

\maketitle

\begin{abstract}
This work develops feasible path trajectories for a coordinated strike with multiple aircraft in a constrained environment.  Using direct orthogonal collocation methods, the two-point boundary value optimal control problem is transcribed into a nonlinear programming problem.  A coordinate transformation is performed on the state variables to leverage the benefits of a simplex discretization of the search domain.  Applying these techniques allows each path constraint to be removed from the feasible search space, eliminating computationally expensive, nonlinear constraint equations and problem specific parameters from the optimal control formulation.  Heuristic search techniques are used to determine a Dubins path solution through the space to seed the optimal control solver.  In the scenario, three aircraft are initiated in separate directions and are required to avoid all constrained regions while simultaneously arriving at the target location, each with a different viewing angle.  A focus of this work is to reduce computation times for optimal control solvers such that real-time solutions can be implemented onboard small unmanned aircraft systems.  Analysis of the problem examines optimal flight paths through simplex corridors, velocity and heading vectors, control vectors of acceleration and heading rate, and objective times for minimum time flight.
\end{abstract}

DISTRIBUTION STATEMENT A: Approved for public release; distribution is unlimited. AFRL-PA-2023-1431

\section*{Nomenclature}

{\renewcommand\arraystretch{1.0}
\noindent\begin{longtable*}{@{}l @{\quad=\quad} l@{}}
$a$             & acceleration \\
$\alpha$        & barycentric coordinate \\
$ct$            & circular orbit time \\
$f$             & state dynamics \\
$ft_n$          & geometric flight time \\
$g$             & constraint function \\
$\gamma$        & change in heading rate \\
$J$             & cost function \\
$p$             & phase \\
$P$             & total number of phases \\
$R$             & turning radius \\
$t_f$           & optimal flight time \\
$td_n$          & flight time difference \\
$\theta$        & heading \\
$\dot{\theta}$  & heading rate \\
$\ddot{\theta}$ & change in heading rate \\
$v$             & velocity \\
$u$             & control \\
$X$             & state vector \\
$x$             & $x$-coordinate position \\
$y$             & $y$-coordinate position \\
$x_c$           & circular orbit center $x$ position \\
$y_c$           & circular orbit center $y$ position 

\end{longtable*}}

\section{Introduction}
\lettrine{M}{ultiple} Small Unmanned Aircraft Systems (SUAS) have been used in coordination to accomplish a wide range of objectives.  Beginning with single aircraft missions, designed for surveillance or target strike, SUAS found a permanent spot in US strategy.  With increased confidence in the trust and capability of unmanned systems, Manned-UnManned Teaming (MUM-T) operations become reality, where a manned aircraft is used in coordination with unmanned vehicle to accomplish a task.  As SUAS capabilities continue to evolve, understanding the coordinated effort on multiple UAS aircraft is required.  Research in the field has examined a multitude of methodologies that can be used for SUAS mission profiles.  The traveling sales problem is well known and has been used for task allocation with multiple SUAS with a goal of collecting information from data pods or simply performing perimeter surveillance around an installation \cite{Humphreys2016a}.  Swarming algorithms have been developed to maximize the probability of accomplishing an objective upon which cooperation amongst the vehicles is not prioritized, or random paths are determined for a single trajectory solution \cite{Torres2017}. And SUAS have been use in MUM-T operations to accomplish a mission with the SUAS while a human interface is maintained in the loop \cite{Hocraffer2017}.  Finally, Xue performed research for flight path trajectory planning, approaching real-time computational speed, while avoiding constraints in  the urban environment \cite{Xue2020}.  For each of these mission profiles, the aircraft must maintain some sense of coordination to assure the assigned task is performed in an efficient manner.  

The work herein develops path solutions for multiple SUAS in constrained environments such that given an initial starting position for all aircraft, a coordinated time of arrival at the objective point may be achieved.  A scenario is developed in which a base perimeter has been established with multiple constraints.  The objective for the team of SUAS is to fly a variable speed, optimal trajectory to the target location, with a simultaneous arrival team, while each UAS approaches the target at a different angle for maximum observation.  Each path is unique to each aircraft and therefore control is given to both the rate of change in the vehicles heading rate as well as the vehicles acceleration.  By achieving these trajectories, a detailed view of the target can be acquired while minimizing observance from the enemy.  Ultimately reducing the time required to make operational decisions on the target region as a more complete vision of the target will be acquired instantaneously, giving a greater impact to accomplishing the mission objectives.

\section{Background}
Path planning algorithms have become more popular as computational efficiencies continue to increase, allowing the gap to be reduced between pre-planned and real-time trajectory development.  As optimal trajectory algorithms have been evolved, research has been conducted for optimal trajectories under a myriad of cost functions, to include, minimizing parameters such as flight time \cite{Marinis2022}, energy \cite{Yacef2017}, noise pollution \cite{Scozzaro2017}, and deviations from a known path due to wind disturbances or sensor noise \cite{khachumov2022}. However, in each of these cases, path constraints in the domain have been kept to a minimum.  These constraints are modeled with either piece-wise or continuous functions with potential for steep gradients, thus drastically increasing computational time for each additional constraint added to the domain.  Previous research evaluated the effects of these constraints on computation time and solution feasibility; illustrating trajectory issues that arise when solving discrete systems \cite{Zollars2017a}.  

Considering the coordination of multiple aircraft, Schuldt and Kurucar examined how teaming of multiple SUAS can maximize a search and rescue effort while significantly reducing the time required to search the space \cite{Schuldt2017}.  Keller's work examined the ability for multi-fixed wing aircraft systems to achieve simultaneous arrival at a target location given unspecified starting positions \cite{Keller2014}.  Specifically, by using geometric path relations through B-splines, simultaneous arrival to a target can be achieved by lengthening the total flight time of each aircraft's trajectory to match that of the longest flight time.  This time differential is then accounted for by varying the velocity profile of each aircraft. 
 However, this work operated in free airspace, without incorporating path constraints in the domain.  

In order to build computational efficient algorithms for constrained gradient based optimization, previous research evaluated simplex discretization techniques employed to remove the physical path constraints from the search space.  By discretizing the space, any simplex that is contained in a physical constraint can be removed from the search space, leaving only unrestricted paths for a gradient based search algorithm to explore.  In order to explore the space efficiently, an A* algorithm is employed to find the Connected Simplex Corridor (CSC), a connected path of simplexes, free of any path constraints, from the starting location to the end location.  Once the CSC is defined, the vehicle dynamics are transcribed into the barycentric coordinate frame to evaluate a continuous, optimal path through the CSC. Kallman's heuristic algorithm developed two-dimensional trajectory solution through a highly constrained environment with computation times in the milliseconds  \cite{Kallmann2010a}.  Leveraging Kallman's work with gradient based optimal control solvers, optimal trajectories were analyzed through 37 constraints with computational time of 2 seconds \cite{Zollars2017b}.  That research was then extended for acceleration control and heading rate of the vehicle, allowing the aircraft to have a smaller turn radius and therefore explore a larger region of the constrained space \cite{Zollars2017}.

To determine the optimal trajectory within the CSC, direct orthogonal collocation techniques were used.  Direct orthogonal collocation methods are used to transcribe the two-point boundary value problem to a Nonlinear Program (NLP) problem \cite{Smith2014, Suplisson2015}.  Gaussian quadrature is used for an exact solution of the state and control at each collocation point with Lagrange polynomial interpolation approximating the state and control between each collocation point \cite{Betts2010}. 

\section{Methodology}
\subsection{Initiating the Algorithm}
The algorithm is initiated by defining the polygonal constraint map.  Each building or path constraint, defined as a keep out region, within the search domain is defined as a polygon.  Additionally, the target region is defined as a fictitious constraint upon which each aircraft will terminate its trajectory.  This target polygon is to have as many sides as there are aircraft in the scenario, with each aircraft being assigned a final trajectory location on an edge of the target polygon, thus allowing for maximum observation of the target location.  With the physical constraints defined, the initial conditions for the aircraft's starting position and initial heading are set along with the final conditions to include the edge of the target polygon, as well as required final heading.  
\subsection{Heuristic Solution}
\begin{figure}[htbp!]
\includegraphics[width=4in]{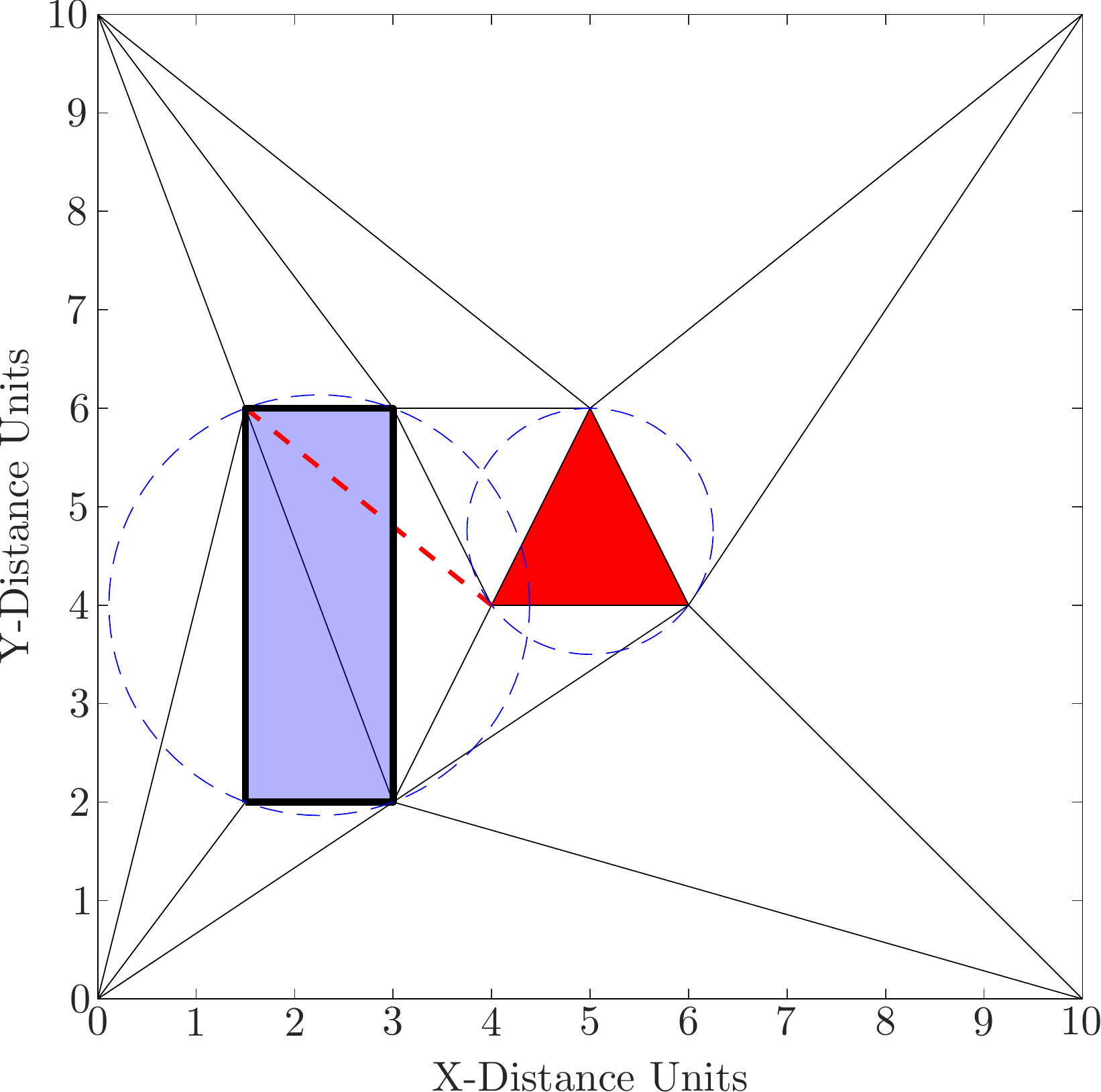}
\centering
\caption{Constrained Delaunay Triangulation}
\label{fig:delaunay}
\end{figure}

With the initial and final position for each aircraft defined, as well as the constraint map, the heuristic solution is solved.  This process is initiated by discretizing the domain into a Constrained Delaunay Triangulation (CDT).  A Delaunay triangulation is defined by Chew as follows \cite{Chew1987}: \\
\textbf{Definition:} Let \textit{S} be a set of points in the plane.  A triangulation \textit{T} is a Delaunay triangulation of \textit{S} if for each edge \textit{e} of \textit{T} there exists a circle with the following properties:
\begin{itemize}
    \item The endpoints of edge \textit{e} are on the boundary of \textit{C}
    \item No other vertex of \textit{S} is in the interior of \textit{C}.
\end{itemize}
Expanding this further for a CDT.  A set of constrained edges is defined as set \textit{G}.  If a vertex point lies within \textit{C} of the constrained triangle from \textit{G}, that vertex point is blocked from one of the vertices of \textit{C}.  Figure \ref{fig:delaunay} illustrates each of these conditions.  First the target simplex is shown as a red triangle, also representing a triangle in the discretization.  A circle is shown to circumvent the three points of the target location, illustrating the basic principles of a delaunay triangulation that no other vertex point can lie within that circle.  The blue shaded region represents a path constraint defined by a set of constrained points \textit{G}.  Through the CDT, it is guaranteed these constrained edges will remain in the simplex set, ultimately allowing for them to be removed from the search space. The circle circumventing the triangle in the constrained region illustrates the CDT differs from the Delaunay triangulation as the vertex point at position $(4,4)$ is clearly inside of this circle, yet the CDC prevents a segment in the discretization from point $(4,4)$ across the constrained edge to point $(1.5,6)$, shown with the thick red dashed segment, thus maintaining the edge of the constrained polygon.  

Following the CDC discretization, a set of simplexes are acquired.  Any simplex that is contained inside a constrained region is removed from the set.  With a feasible search space defined, an A* algorithm is performed to find a connected set of simplexes from the simplex that contains the starting position, to the defined edge of the target simplex, also described as a CSC.  The A* search algorithm gets its roots from Dynamic Programming and is a subset of Dijkstra's search algorithm \cite{Dechter1985}.  The difference in that the A* algorithm can be terminated when an adequate solution has been discovered, even if it may not be optimal, whereas Dijkstra's algorithm searches the entire domain, ensuring the optimal solution is found.  The A* cost function is defined as
\begin{equation}
f(n) = g(n)+h(n)
\end{equation}
where $g(n)$ is the cost evaluated at the current node, and $h(n)$ is the defined as the cost-to-go.  In this work, $g(n)$ is defined as the euclidean distance from the current location to each of the vertex points of the adjoining simplexes. The cost to go, $h(n)$, is defined as the euclidean distance from one of the adjoining vertex points to the final location.  In the case where a segment, defined in the cost-to-go by the euclidean distance, passes through a constrained edge, the cost-to-go is calculated by including each edge of the constrained simplex until a free path to the final location can be achieved. The main disadvantage of the A* algorithm is its dependency on the heuristic chosen for $h(n)$, as this heuristic can change the path solution acquired.

\begin{figure}[htbp!]
\includegraphics[width=3in]{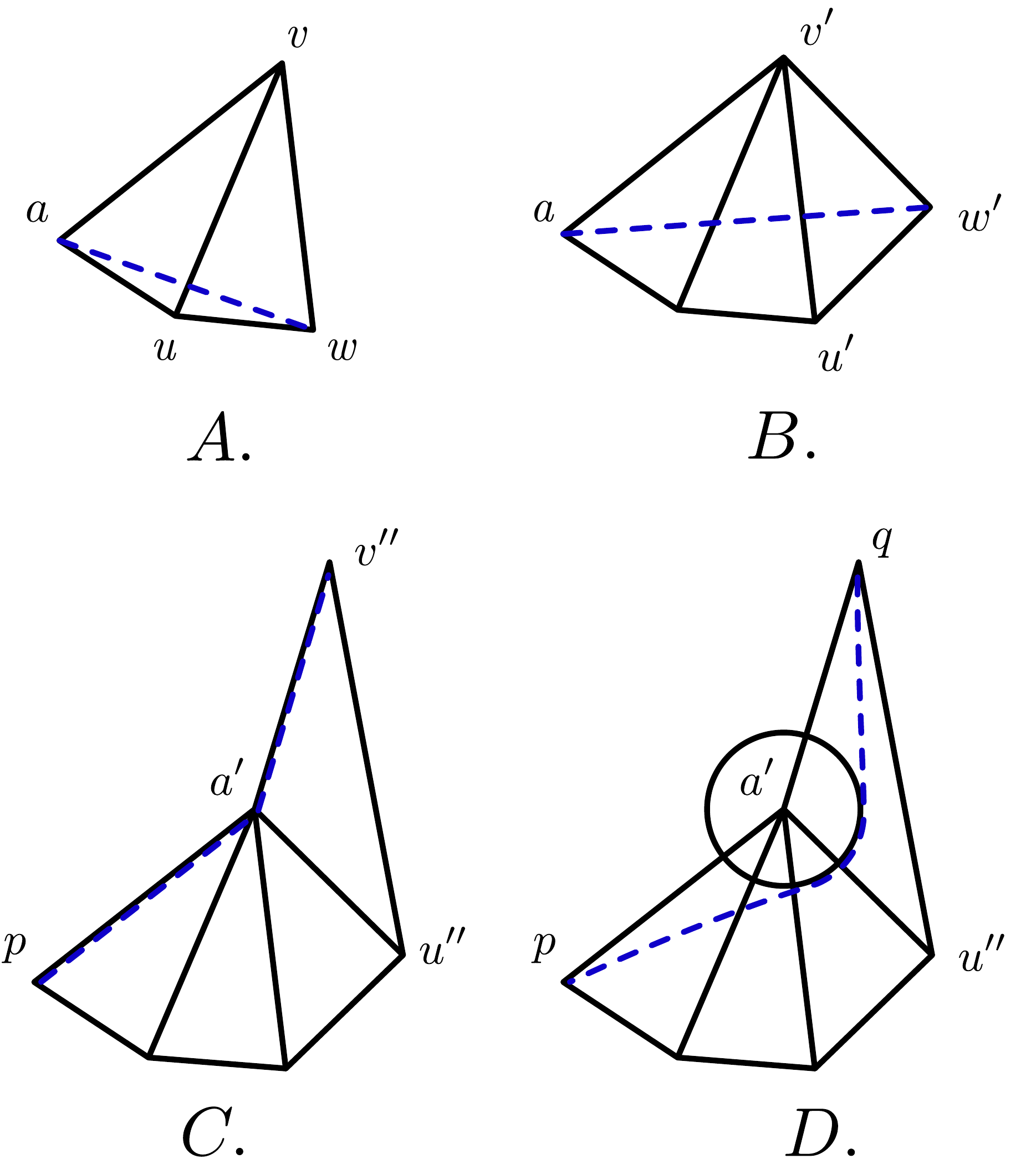}
\centering
\caption{Funnel algorithm}
\label{fig:funnel}
\end{figure}

Once a CSC is obtained for each aircraft, a Dubins path solution is determined based on Lee and Preparata \cite{Lee1984}, and Chazelle's \cite{Chazelle1982} funnel algorithm, as cited by Hershberger \cite{Hershberger1994}.  The algorithm is initiated by defining an initial and final location through the CSC as, $[p,q]$.  A series of straight line paths and minimum radius curves are connected to given an overall path solution.  This solution provides a feasible geometric path through the constrained region, but does not provide any information on state aircraft dynamics or rate limited control inputs.  For further insight into these parameters, gradient-based optimal control algorithms are derived.  Figure \ref{fig:funnel}, depicts this process starting at a vertex point $a$.  

With the CSC defined, Figure \ref{fig:funnel}A shows vertex points $v$ and $u$  define a shared edge between the two simplexes, while $w$ defines the third vertex of the next simplex in the CSC.  If a straight line path from $a$ to $w$ is feasible, that path is temporarily chosen as the current Dubins path solution.  Evaluating the third simplex in figure \ref{fig:funnel}B, the CSC once again evaluates the unshared vertex point, defined here as $w'$.  If the straight line path is feasible, the Dubins path solution is updated.  If the straight line path is not feasible, as shown in Figure \ref{fig:funnel}C, the closest shared vertex is chosen as the new apex point, now defined as $a'$ and the funnel algorithm continues to the next simplex in the CSC.  In order for the Dubins path to be feasible, a circle is placed over each apex point in the path defined by the minimum turning radius, R, of the vehicle,
\begin{equation}
R = \frac{v_\mathrm{min}}{\dot{\theta}_\mathrm{min}}
\end{equation}
where $v_\mathrm{min}$ defines the minimum aircraft velocity and $\dot{\theta}_\mathrm{min}$ defines the aircraft's minimum turn rate. Tangent points to this curve are calculated for each apex location and the Dubins path solution is determined as shown in \ref{fig:funnel}D.

\subsection{Vehicle Dynamics}
 This problem is solved in the optimal control general purpose solver, GPOPS II, with an sequential quadratic programming NLP solver.  The objective of this work is to acquire a feasible path through a constrained environment in minimum time with a focus on reducing computation time for real-time optimal control of SUAS.  

The two-dimensional aircraft state dynamics are given as 
\begin{equation}
\label{eq:dyn1}
 \mathbf{f}(\mathbf{x}(t), \mathbf{u}(t), t) = \begin{bmatrix} \dot{x}(t) \\ \dot{y}(t) \\ \dot{\theta}(t) \\ \ddot{\theta}(t) \\ \dot{v}(t) \end{bmatrix} =
 \begin{bmatrix}
 v \cos(\theta(t)) \\
 v \sin(\theta(t)) \\
 \dot{\theta}(t) \\
\gamma(t) \\
 a(t)
 \end{bmatrix}.
\end{equation}
Control for the vehicle is the change in heading rate and acceleration,
\begin{equation}
\label{eq:cont}
\mathbf{u}(t) = [\gamma(t), a(t)].
\end{equation}
These dynamics are translated to barycentric coordinates for propagation of the state and control through a discretized simplex mesh.
The simplex mesh is constructed with a CDT formed from the polygonal constraints in the domain \cite{Zollars2017a}.  Given an initial and final state for position, A* search algorithm is applied to find a set of simplexes in the space that assemble a feasible path solution known as a CSC.  Once the series of triangles in the CSC are known, a Dubins path solution is obtained to seed the optimal control solver.  

In order for each aircraft to converge on the final target at the same time in an optimal fashion, an orbit of the starting location may be performed for each aircraft following the aircraft with the longest path.  Given the Dubins path result for each aircraft, the longest path is determined and the flight time difference between longest path and each subsequent aircraft is calculated.  Given this flight time difference for each subsequent aircraft, the number of potential orbits of the starting location is determined.  For these orbits, the aircraft turn radius is fixed, restricting each vehicle to a flight path of a fixed radius circle.  The time difference between the aircraft with the longest path and all other aircraft, is then consumed in this initial orbit, thus allowing for a minimum time flight profile from the end of the orbit to the target location for each aircraft.  

\subsection{Optimal Control Problem}
The Dubins path solution provided the $(x,y)$ position and max velocity time of flight for each aircraft.  Direct orthogonal collection methods are now used to determine optimal flight paths, providing full state trajectories of the solution.  To acquire efficient flight path trajectories for each aircraft of the simulation, two optimal control problems are solved. The first optimal control problem evaluates the time difference between the longest flight path and each successive flight path to determine if an aircraft is required to orbit the starting location before beginning its path to the target.  By evaluating the flight time produced from the Dubins path solution for each vehicle, the minimum flight path time for each vehicle is determined.  The vehicle with the longest flight time, designated as aircraft 1, will bypass the orbit stage and immediately transition to the target.  The remaining aircraft will either begin a circular path around the starting location or begin to transition to the target based on the time difference between the their Dubins path solution and the flight time of aircraft 1.  This logic is determined by the flight times of the Dubins path solution, defined as
\begin{equation}
\begin{aligned}
ft_1 &= Flight \: time \: of \: longest \: path \\
ft_2 &= Flight \: time \: of \: second \: longest \: path \\
ft_N &= Flight \: time \: of \: shortest \: path.
\end{aligned}
\end{equation}
For aircraft's $2$ through $N$, the time difference, $td_n, \: \forall \: n \in [2 \dots N]$, is calculated and number of circular orbits around the starting location is determined based on the minimum and maximum time to complete one constant radius orbit given specifications of the vehicles velocity.  
\begin{equation}
\begin{aligned}
1 \: \text{orbit} &:= td_n \leq ct_\mathrm{min} \\
2 \: \text{orbits} &:= 2(ct_\mathrm{max}) \leq td_n \leq 2(ct_\mathrm{min}) \\
N \: \text{orbits} &:= N(ct_\mathrm{max}) \leq td_n \leq N(ct_\mathrm{min})
\end{aligned}
\end{equation}
where $ct_\mathrm{min}$ and $ct_\mathrm{max}$ is the time required to fly a single circular path of fixed radius given the maximum or minimum velocity of the aircraft respectively.  Ultimately, the objective is to minimize transition flight time between the starting point and the target for all aircraft. 

The first optimal control problem has a starting and final location for each vehicle at $(x_0,y_0)$, and an initial and final heading of $\theta_0$.  The dynamics are given in Equation (\ref{eq:dyn1}) with control defined in Equation (\ref{eq:cont}).  In order to maintain a circular orbit, a path constraint is added to maintain a radius from the center point location, $(x_c,y_c)$, of the orbit equal to the vehicle's maximum velocity minimum turning radius, $r$, defined as
\begin{equation}
\mathbf{g}(\mathbf{x}(t),\mathbf{u}(t),t) = (x(t)-x_c)^2+(y(t)-y_c)^2-r^2.
\end{equation}
Finally, the cost function was designed to minimize the acceleration control on the vehicle in order to maintain efficient fuel consumption and provide a smooth heading rate vector, defined as
\begin{equation}
J_1= \int{a^2(t)}dt.
\end{equation}

\begin{figure}[htbp!]
\includegraphics[width=3.5in]{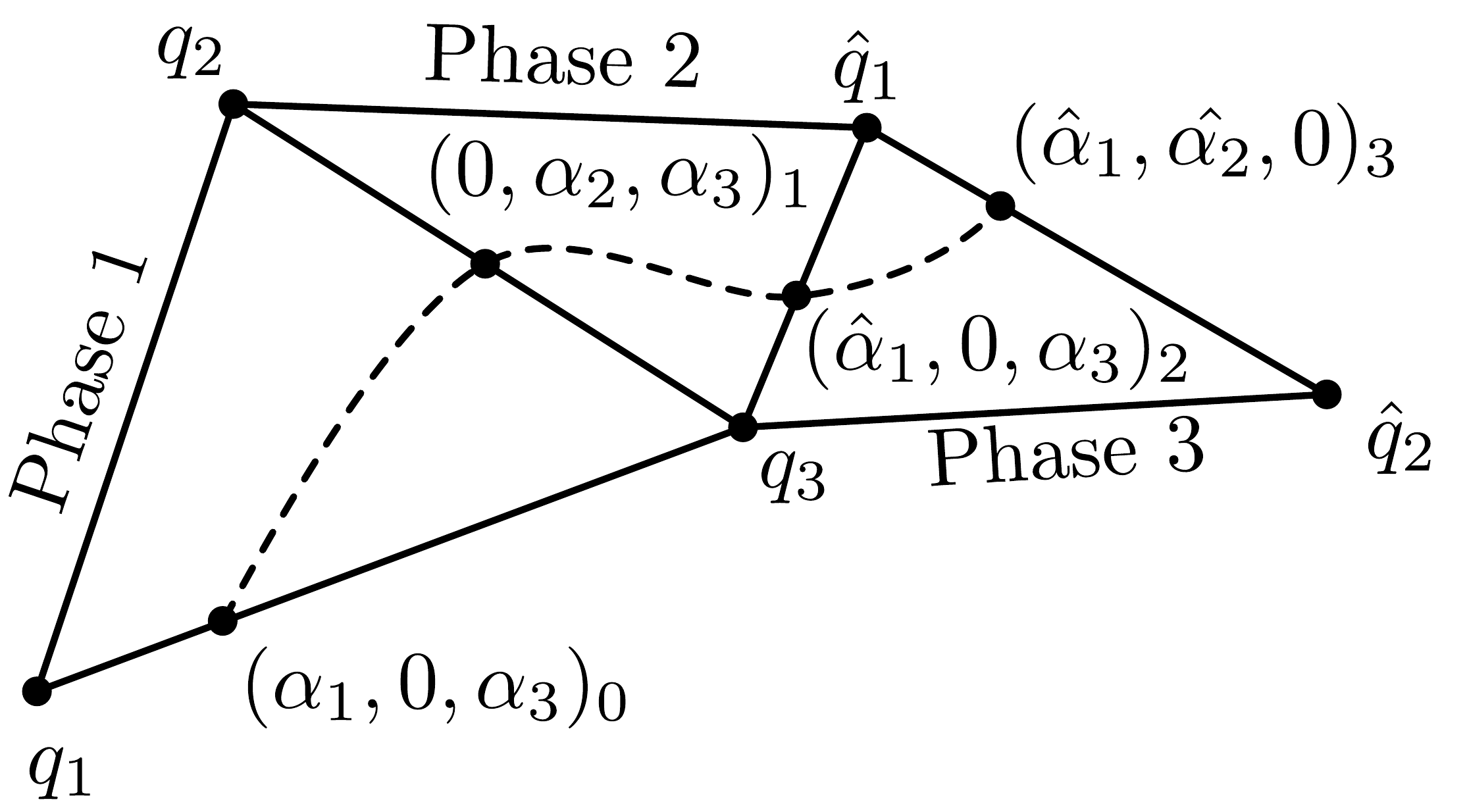}
\centering
\caption{Phase solution through barycentric coordinate frame}
\label{fig:phases}
\end{figure}

The second optimal control problem is initiated by first establishing the start location and the target simplex.  The target simplex is defined as a polygon of sides equal to the number of aircraft in the mission.  Each aircraft will end its trajectory on a predefined edge of the target simplex.  The states, defined in Equation (\ref{eq:dyn1}), are translated to the barycentric coordinate frame.  This translation allows for propagation of the state trajectory through each individual simplex, where the vehicle position is a coordinate triplet, defined by the weight to each vertex point of the active simplex \cite{Zollars2017a}.  With the CSC determined previously with the A* algorithm, each simplex of the CSC is optimally connected through event constraints within a phased structure of the General Purpose Optimal Pseudospectral Solver (GPOPS-II).  The state vector represents the position in barycentric coordinates, heading angle, heading angle rate, and velocity, defined by the dynamic equations 
\begin{align}
\dot{\alpha}_1^{(p)}(t) &=\frac{ (y_2 - y_3)\cos(\theta(t))  + (x_3-x_2)\sin(\theta(t))}{\det(T)}v(t)  \label{eq:DEi} \\
\dot{\alpha}_2^{(p)}(t) &=\frac{ (y_3 - y_1)\cos(\theta(t))+ (x_1-x_3)\sin(\theta(t))}{\det(T)} v(t)  \\
\dot{\alpha}_3^{(p)}(t) &= -\dot{\alpha}_1^{(p)}(t)-\dot{\alpha}_2^{(p)}(t) \\
\dot{\theta}^{(p)}(t) &=\theta(t) \\  
\ddot{\theta}^{(p)}(t) &= \gamma(t) \\
\dot{v}^{(p)}(t) &= a(t)  \label{eq:DEf}
\end{align}
$\forall p \in [1...P-1]$, where the superscript, $p$, defines the phase of the solution, also referred to as the active simplex within the CSC while $P$ defines the total number of phases.  This process is illustrated in Figure \ref{fig:phases}.  

The resulting state vector becomes
\begin{equation}
\label{eq:state}
    \mathbf{X}(t) = [\alpha_1, \alpha_2, \alpha_3, \theta, \dot{\theta},v]^\top.
\end{equation}
As each phase of the optimal trajectory is solved, event constraints are implemented to assure continuity of the state, control and time, defined as 
\begin{align}
\mathbf{X}_0^{(p+1)}-\mathbf{X}_f^{(p-1)}&=0 \: \: \forall p\in [2 \dots P] \\
\mathbf{u}_0^{(p+1)}-\mathbf{u}_f^{(p-1)}&=0 \: \: \forall p\in [2 \dots P] \\
t_0^{(p+1)}-t_f^{(p-1)}&=0 \: \: \forall p\in [2 \dots P].
\end{align}

The control for the SUAS is on the change in heading rate and acceleration as defined previously, but now defined in the phased structure as
\begin{equation}
\mathbf{u}^{(p)}(t) = (\gamma^{(p)}(t), a^{(p)}(t)).
\end{equation}
Bounds are applied on the states, control, and time to limit the search domain of the NLP solver, enforced as 
\begin{align}
\label{eq:BDSI}
0 &\leq \alpha_1^{(p)} \leq 1 \\  
0 &\leq \alpha_2^{(p)} \leq 1 \\
0 &\leq \alpha_3^{(p)} \leq 1 \\
 | \theta^{(p)} | &\leq \pi \\
 | \dot{\theta}^{(p)} | &\leq 45 \enspace \text{deg/s}  \\
 10 \enspace \text{ft/s} &\leq v^{(p)}  \leq 30 \enspace \text{ft/s} \\
 | \gamma^{(p)} | &\leq 10  \enspace \text{deg/s}^2  \\
 | a^{(p)} | &\leq 5  \enspace \text{deg/s}^2  \\
 0 \leq t^{(p)} &\leq \frac{\text{edge}_\mathrm{max}^{(p)}}{v} \label{eq:BDSf}
\end{align}
$\forall p \in [1...P]$, where $\text{edge}_\mathrm{max}$ represents the longest edge in the defined simplex.

The objective of the second optimal control problem is to solve the flight trajectories for all vehicles with a coordinated arrival time, over all phases subject to the dynamic constraints defined in Equations \ref{eq:DEi}-\ref{eq:DEf}, and parameter bounds of Equations \ref{eq:BDSI}-\ref{eq:BDSf}.  The cost function is defined as
\begin{equation}
J_2 = {|t_{f_1}-t_{f_n}|}
\end{equation}
where $t_{f_1}$ is defined as the summation of the flight time of each simplex within aircraft 1's CSC, as
\begin{equation}
t_{f_1}=\int_{t_{0}^{(p)}}^{t_f^{(p)}}dt \: \: \forall p\in[1...P].
\end{equation}
For the remaining aircraft, $t_{f_n}$ is defined as the summation of the initial orbit flight time summed with the flight time of each simplex within that aircraft's CSC, as
\begin{equation}
t_{f_n}=J_{1_n}+\int_{t_{0}^{(p)}}^{t_f^{(p)}}dt \: \: \forall p\in[1...P].
\end{equation}
This cost function is designed to minimize the time difference between the first aircraft's optimal trajectory to each successive aircraft's total flight time, thus providing a solution upon which all aircraft arrive at the target location simultaneously. 

\section{Multi-UAS Scenario}
\begin{figure}[htbp!]
\includegraphics[width=4in]{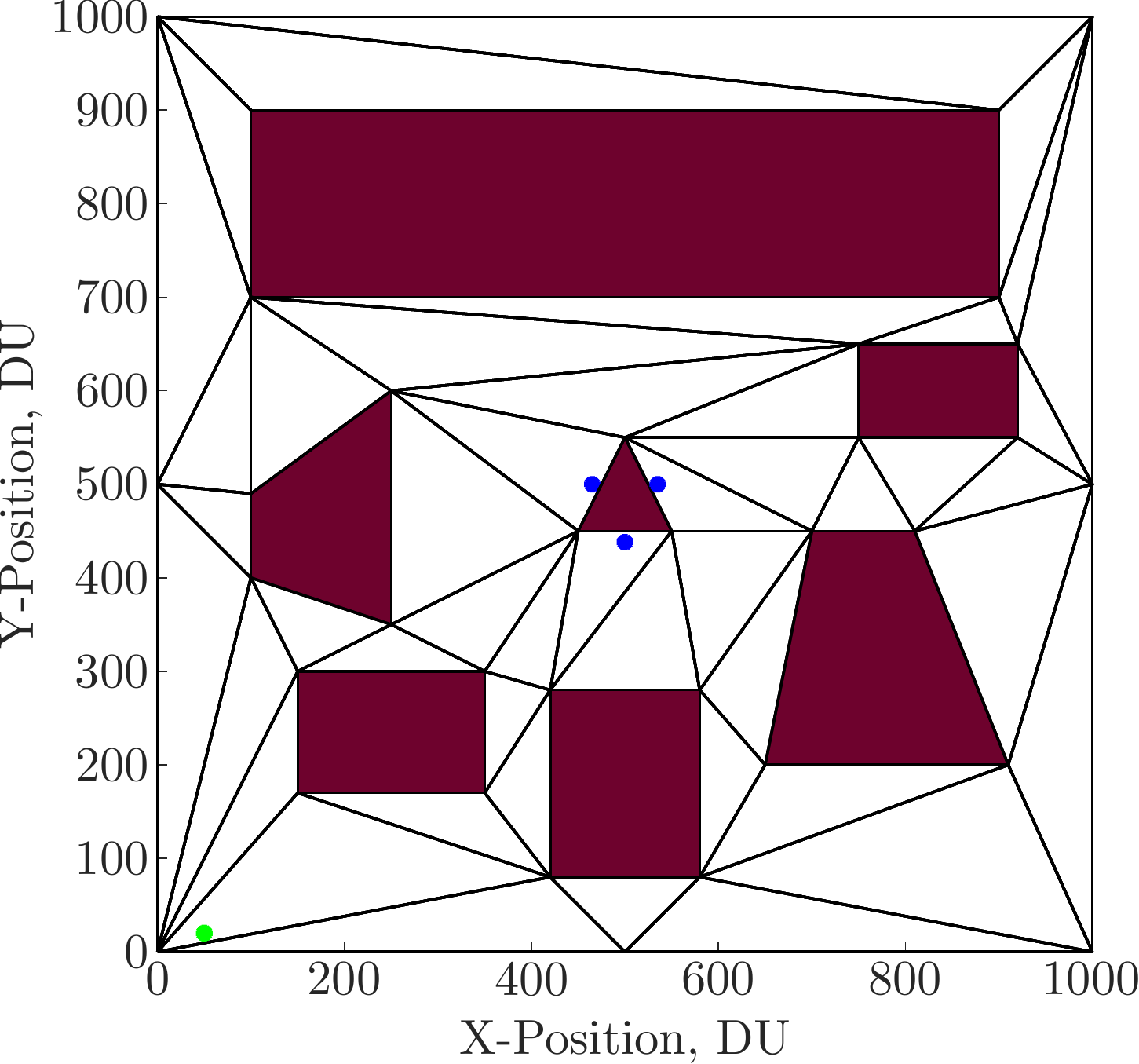}
\centering
\caption{Simplex discretization}
\label{fig:simplex}
\end{figure}

The scenario is constructed for a target surrounded by multiple path constraints.  Each path constraint resembles a physical structure, or keep-out region.  In order for the trajectory of the aircraft to maintain a safe distance from the building structure throughout the trajectory, a building offset is included in the polygonal constraint, therefore allowing a trajectory to safely approach coincidence to the boundary of the simulated keep-out region.  Three identical aircraft are initiated at the same location, each with a zero heading angle and flying at max velocity.  It is assumed the aircraft are separated in altitude for airspace deconfliction.  The target is located at the center of the field and is encompassed by a three-sided simplex.  Each aircraft is assigned a terminal position and heading corresponding to an individual side of the target simplex.  The initial constraint map is detailed in Figure \ref{fig:simplex}.
The starting position of each aircraft is indicated by the green asterisk at location $(50,20)$, with the final location at the blue asterisks surrounding the center simplex.  Each dark red simplex represents a keep-out region.  The input to the algorithm requires the initial and final position and heading of each aircraft, the target location, coordinates for each vertex of the keep-out regions, and the bounds on the domain.  The initial step is for a CDT to be performed on the domain.  This triangulated discretization is shown in Figure \ref{fig:simplex}.

\begin{figure}[htbp!]
\includegraphics[width=4in]{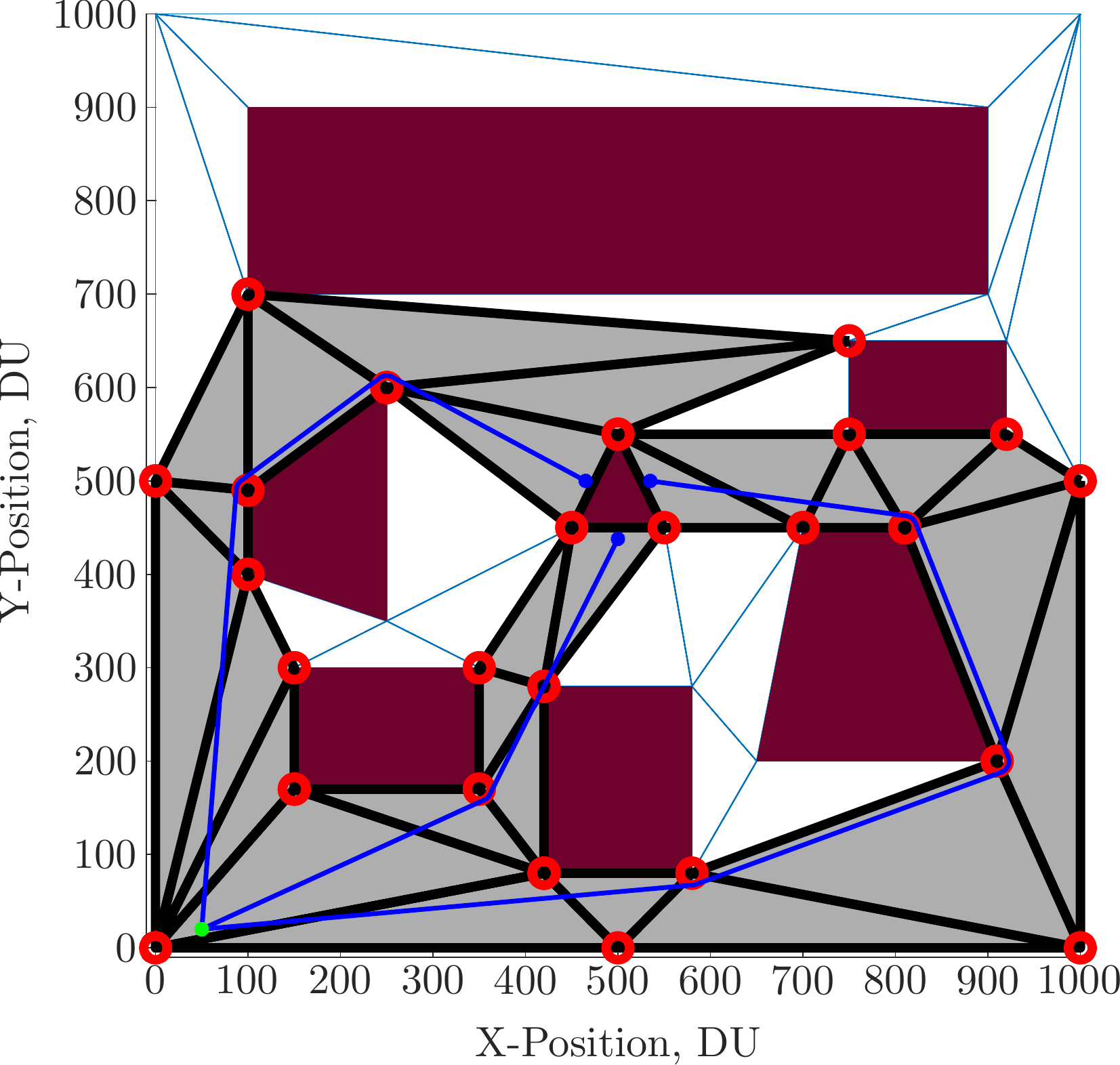}
\centering
\caption{Multi-aircraft Dubins path solution}
\label{fig:dubins}
\end{figure}

With the simplex map defined, an A* search algorithm is executed for the first aircraft to find a feasible set of simplexes from the initial position to the desired edge of the target simplex, also known as the CSC.  For the remaining aircraft, a modified A* search algorithm is performed such that the desired CSC is produced for each remaining aircraft in the scenario.  Each of the CSC's are shown in Figure \ref{fig:dubins}, resulting in each aircraft progresses to the target in different directions, minimizing the threat of detection and maximizing the viewing angle of the target.   
Both the Dubins path solution and the minimum turn radius circular offsets can be seen as the solid blue lines and red circles on the edge of each constraint, respectively.  As the Dubins path is calculated via the funnel algorithm, if the straight line path intersects a constrained edge, the minimum turning radius circle is utilized.  However, it should be noted the minimum turning radius circles are only required if a constrained edge is violated. 
 A straight line path that only crosses through the minimum turning radius circle remains a feasible path due to the building offset included in each physical constraint, as seen at location $[420,280]$.  This geometric path solution represents a feasible flight path and is used to seed the optimal control solution. Flight times are determined assuming each vehicle flies this geometric path at maximum speed and are given in Table \ref{table:dubinsTimes}.
\begin{table}[htbp!]
\centering
\caption{Dubins Path Flight Times}
\begin{tabular}{|| c | c ||} 
 \hline
 \textbf{Vehicle} &  \textbf{Time (sec)} \\ 
 \hline
 Aircraft 1 & $51.53$ \\
 \hline
 Aircraft 2 & $32.41$ \\
 \hline
 Aircraft 3 & $22.34$ \\
 \hline
\end{tabular}
\label{table:dubinsTimes}
\end{table}

Previous work under this methodology has shown the Dubins path solution to be an upper bound on the flight path time, as the optimal solution is permitted to intersect the minimum radius turn circles \cite{Zollars2017c}.  Therefore the time difference between the max flight time and all other trajectories is calculated and the time difference for each aircraft is consumed flying a fixed radius, constant altitude orbit at the start location.  In order to efficiently consume the entire time difference, control is given to the vehicles acceleration and the change in heading rate as described in Equation \ref{eq:cont}.  This first optimal control problem is defined as a fixed initial state, fixed final final state, fixed final time, with the objective to minimize control on acceleration.  To maintain consistency with the Dubins path solution,  the same boundary conditions and rate limits on the state parameters are applied to the initial and final state and the control of the aircraft.  For the circular orbit, the aircraft's initial and final conditions are
\begin{eqnarray}
&[x_0, y_0,\theta_0,v_0] = [20, 50, 0, v_\mathrm{max}] \\
&[x_f, y_f,\theta_f,v_f] = [20, 50, 0, v_\mathrm{max}] 
\end{eqnarray}

\begin{figure}[htbp!]
\includegraphics[width=4in]{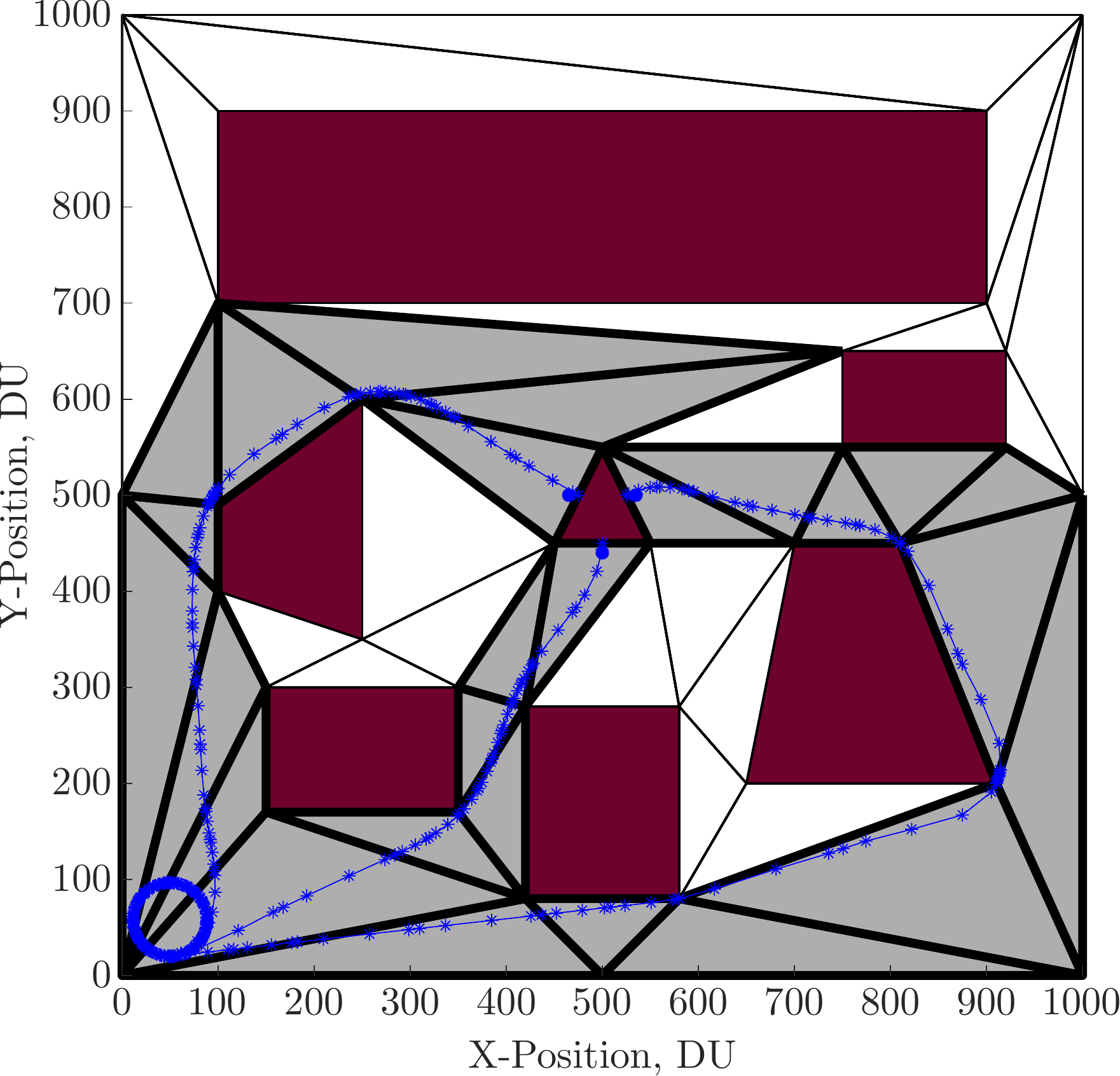}
\centering
\caption{Multi-aircraft optimal path solution}
\label{fig:optimal}
\end{figure}

The second optimal control problem solves the phased optimal trajectory given each aircraft's CSC with the objective of simultaneous arrival at the target location. This path solution is shown in Figure \ref{fig:optimal}.
Here each collocated solution is shown along with the initial orbit for aircraft $2$ and $3$. The far right path defines the longest trajectory and shows the minimum time flight for the CSC.  The flight time allocated for the middle and left paths are dictated by the Dubins path solution and therefore the trajectory solution solved is not a minimum time path, rather a path solution designed to assure a coordinated, simultaneous arrival at the target location. The optimal flight times are shown in Table \ref{table:optTimes}.
\begin{table}[htbp!]
\centering
\caption{Optimal Path Flight Times}
\begin{tabular}{|| c | c | c ||} 
 \hline
 \textbf{Vehicle} &  \textbf{Orbit Time (sec)} &  \textbf{Transit Time (sec)} \\ 
 \hline
 Aircraft 1 & $0$ & $48.72$\\
 \hline
 Aircraft 2 & $16.31$ & $32.41$\\
 \hline
 Aircraft 3 & $26.37$ & $22.35$\\
 \hline
\end{tabular}
\label{table:optTimes}
\end{table}
Aircraft $1$ had the longest trajectory and therefore had zero orbit time. Aircraft $2$ and $3$ each performed two circular orbits before transiting to the target location.  Computationally, the algorithm solves in $16.71$ seconds on an iMac, Apple M1 chip, 8 GB memory computer. Although minimizing computational speed is important for this work, it was not an objective of this work.  For a pre-planning flight trajectory algorithm, the computation speed for algorithm is within the bounds for this technology development.  

\begin{figure}[htbp!]
\includegraphics[width=4.5in]{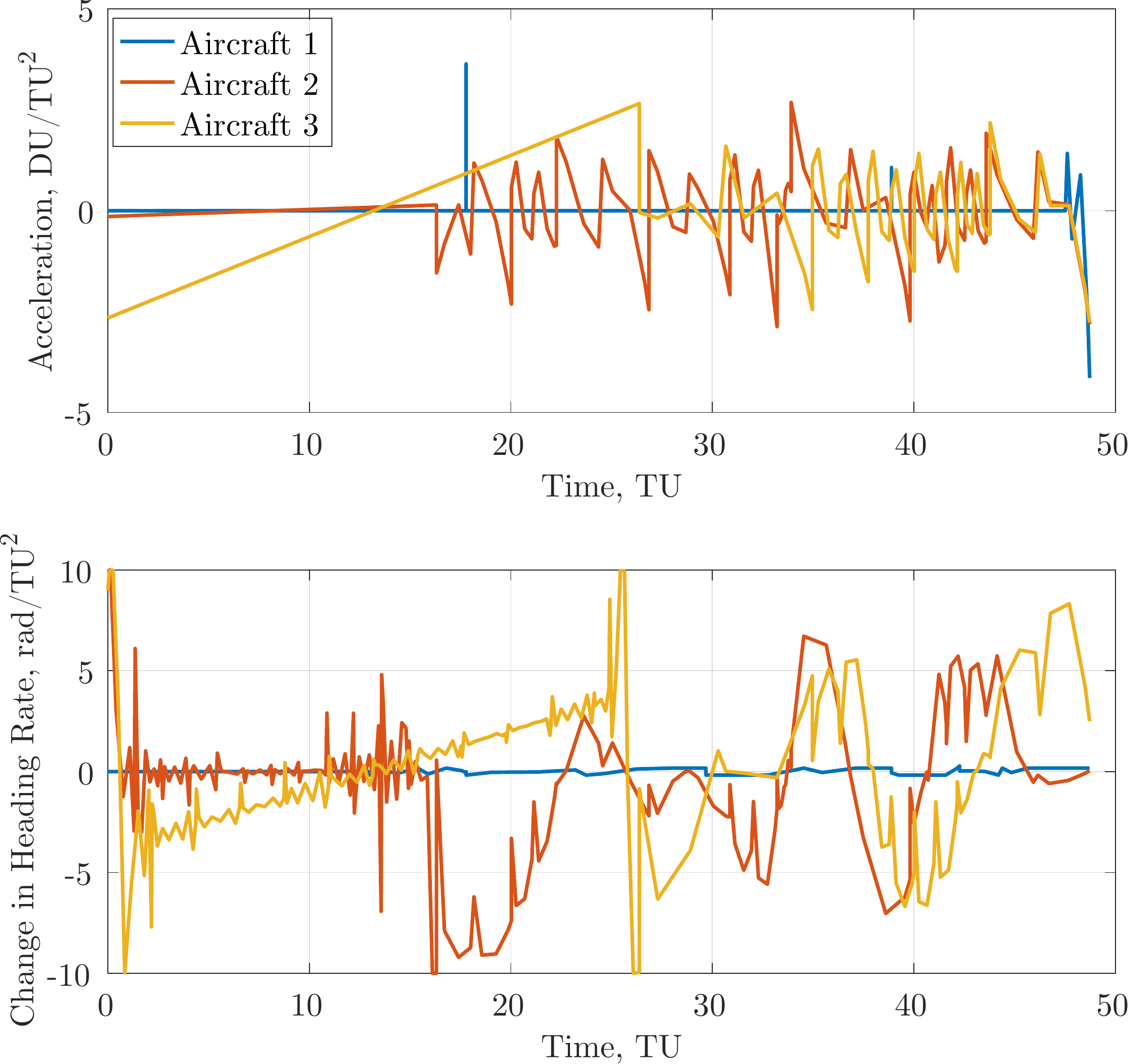}
\centering
\caption{Optimal control profile}
\label{fig:control}
\end{figure}

\begin{figure}[htbp!]
\includegraphics[height=3.5in]{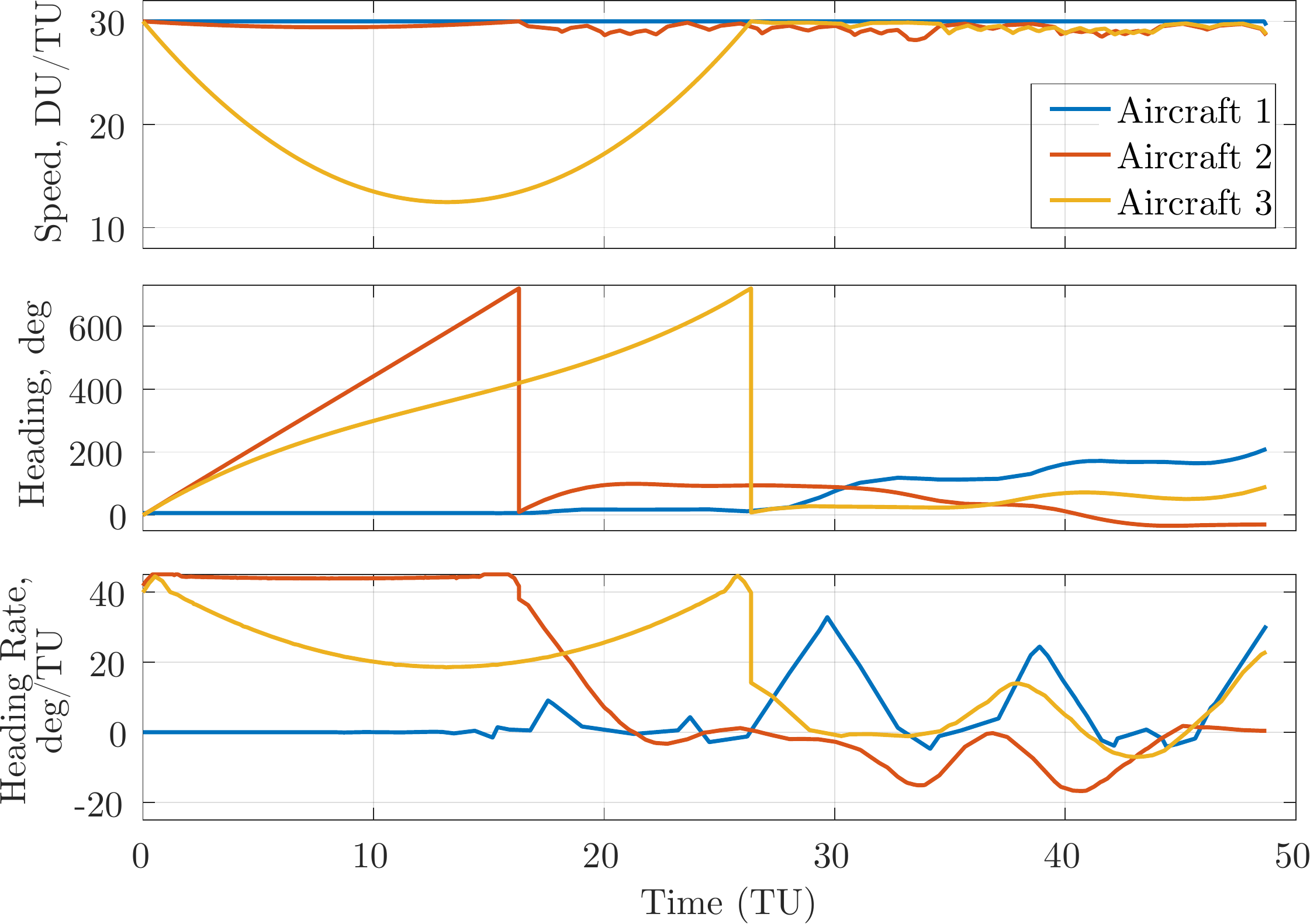}
\centering
\caption{Optimal state profiles}
\label{fig:states}
\end{figure}

The control calculated for the optimal solution of each of the three aircraft is shown in Figure \ref{fig:control}.
Although this is the control calculated for the optimal solution, a rate limited, smooth curve is desired.  For application, the input to the aircraft would be the velocity and heading rate vector shown in Figure \ref{fig:states}.
Here, the entire flight trajectory for each aircraft can be seen.  The far right trajectory, aircraft $1$, had the longest trajectory as it maintained a maximum velocity throughout the flight path.  Aircraft $2$ represents the far left trajectory and it can be seen to fly close to a max velocity circular orbit before transitioning to the target location.  The two circular orbits can be seen in the heading angle and max heading rate, culminating at the $16$ second mark.  Aircraft $3$, represented in the center trajectory, had a time difference of $26.37$ seconds, of which are accounted for in two variable speed circular orbits.  For this simulation, an aircraft could accomplish a circular orbit in $8$ seconds. Therefore, aircraft $3$ could have flown $3$ circular orbits close to maximum speed, vice $2$ slower orbits.  For energy saving purposes, this control algorithm defaulted to a minimum number of circular orbits for each aircraft.  The variable speed of the circular orbit for aircraft $3$ is clearly shown in each of the velocity, heading, and heading rate trajectories.

\section{Analysis \& Conclusion}

The strength of this technique resides in the mesh discretization  removing the polygonal constraints from the search space.  Operationally, this path planning algorithm would be used for pre-planning multi-flight trajectories in which reducing computation speed is a necessity, but convergence times for real-time operations are not required. Results show optimal state and control trajectories for each aircraft, culminating with a coordinated, simultaneous arrival of the target location.  Transformation to barycentric coordinates allow for a phased approach in which all keep-out regions are eliminated from the search space, leaving an unconstrained search domain through a simplex channel, significantly reducing computational speed for highly constrained fields.  

Using direct orthogonal collocation methods for optimal control, within the software package GPOPS-II, multiple optimal aircraft trajectories can now be determined.  To seed the optimal control algorithm, a heuristic approach is taken using an A* search algorithm to determine a feasible set of simplexes form the initial point to the target the location.  With the CSC defined, a geometric Dubins path solution is calculated.  By evaluating the differences in all the flight times of each individual flight path, the first optimal control problem is solved for each aircraft except for that of the longest trajectory.  The aircraft with the longest path calculates its optimal trajectory to the target solution.  The remaining aircraft consume the flight time difference in an orbit at the starting location.  Once an aircraft has consumed all of its time difference, the optimal path solution to the target location is calculated such that the time arrival at the target is coincident for all vehicles.

The utility of this work is essential for multi-aircraft path planning around keep out regions that could be defined as buildings, encampments, or regions where noise must be reduced or eliminated.  This optimal control algorithm allows for flexibility in coordinating aircraft to get the maximum exposure to a target while minimizing flight time through a constrained region.  Future work with simplex discretization will continue to focus on multi-aircraft operations.  Consideration will be given to cost functions designed to minimize energy on a vehicle as well as progressing through simplex channels with a formation of SUAS.  

{\scriptsize
\bibliography{IEEEabrv,Zollars6Feb}}

\begin{thebibliography}{23}
\newcommand{\enquote}[1]{``#1''}
\providecommand{\natexlab}[1]{#1}
\providecommand{\url}[1]{\texttt{#1}}
\providecommand{\urlprefix}{URL }
\expandafter\ifx\csname urlstyle\endcsname\relax
  \providecommand{\doi}[1]{\discretionary{}{}{}https://doi.org/#1}\else
  \providecommand{\doi}[1]{\discretionary{}{}{}\urlstyle{rm}\url{https://doi.org/#1}}\fi

\bibitem[{Humphreys and Cobb(2016)}]{Humphreys2016a}
Humphreys, C.~J., and Cobb, R.~G., \enquote{A Hybrid Optimization Technique
  Applied to the Intermediate-Target Optimal Control Problem,} \emph{Global
  Journal of Technology and Optimization}, Vol.~7, 2016.
\newblock \doi{10.4172/2229-8711.1000200}.

\bibitem[{Torres and Dehn(2017)}]{Torres2017}
Torres, S., and Dehn, J., \enquote{Wind optimal trajectories for UAS and light
  aircraft,} \emph{2017 IEEE/AIAA 36th Digital Avionics Systems Conference
  (DASC)}, 2017, pp. 1--8.
\newblock \doi{10.1109/DASC.2017.8102076}.

\bibitem[{Hocraffer and Nam(2017)}]{Hocraffer2017}
Hocraffer, A., and Nam, C.~S., \enquote{A meta-analysis of human-system
  interfaces in unmanned aerial vehicle (UAV) swarm management,} \emph{Applied
  Ergonomics}, Vol.~58, 2017, pp. 66--80.
\newblock \doi{10.1016/j.apergo.2016.05.011}.

\bibitem[{Xue and Wei(2020)}]{Xue2020}
Xue, M., and Wei, M., \enquote{Small {UAV} flight planning in urban
  environments,} \emph{AIAA AVIATION 2020 FORUM}, Vol. 1 PartF, American
  Institute of Aeronautics and Astronautics Inc, AIAA, virtual event, 2020, pp.
  1--12.
\newblock \doi{10.2514/6.2020-2890}.

\bibitem[{Marinis et~al.(2022)Marinis, Iavernaro, and Mazzia}]{Marinis2022}
Marinis, A.~D., Iavernaro, F., and Mazzia, F., \enquote{A minimum-time
  obstacle-avoidance path planning algorithm for unmanned aerial vehicles,}
  \emph{Numerical Algorithms}, Vol.~89, 2022, pp. 1639--1661.
\newblock \doi{10.1007/s11075-021-01167-w}.

\bibitem[{Yacef et~al.(2017)Yacef, Rizoug, Bouhali, and Hamerlain}]{Yacef2017}
Yacef, F., Rizoug, N., Bouhali, O., and Hamerlain, M., \enquote{Optimization of
  Energy Consumption for Quadrotor UAV,} \emph{International Micro Air Vehicle
  Conference and Flight Competition (IMAV)}, 2017, pp. 215--222.
\newblock \urlprefix\url{http://www.imavs.org/pdf/imav.2017.31}.

\bibitem[{Scozzaro et~al.(2019)Scozzaro, Delahaye, and Vela}]{Scozzaro2017}
Scozzaro, G., Delahaye, D., and Vela, A.~E., \enquote{Noise Abatement
  Trajectories for a UAV Delivery Fleet,} \emph{9th SESAR Innovation Days},
  2019, pp. 1--9.
\newblock \urlprefix\url{https://hal-enac.archives-ouvertes.fr/hal-02388280v3}.

\bibitem[{Khachumov and Khachumov(2022)}]{khachumov2022}
Khachumov, M., and Khachumov, V., \enquote{Optimization Models of UAV Route
  Planning For Forest Fire Monitoring,} \emph{Proceedings - 2022 International
  Russian Automation Conference, RusAutoCon 2022}, Institute of Electrical and
  Electronics Engineers Inc., 2022, pp. 272--277.
\newblock \doi{10.1109/RusAutoCon54946.2022.9896260}.

\bibitem[{Zollars and Cobb(2017)}]{Zollars2017a}
Zollars, M.~D., and Cobb, R.~G., \enquote{Simplex Methods for Optimal Control
  of Unmanned Aircraft Flight Trajectories,} \emph{ASME Dynamics Systems and
  Controls Conference}, 2017, pp. 1--10.
\newblock \doi{10.1115/DSCC2017-5031}.

\bibitem[{Schuldt et~al.(2017)Schuldt, Kurucar, and Schuldt}]{Schuldt2017}
Schuldt, D.~W., Kurucar, J., and Schuldt, D., \enquote{Maritime Search and
  Rescue via Multiple Coordinated UAS,} Tech. rep., MIT Lincoln Laboratory, 6
  2017.
\newblock \urlprefix\url{https://apps.dtic.mil/sti/citations/AD1030377}.

\bibitem[{Keller et~al.(2014)Keller, Thakur, Dobrokhodov, Jones, Likhachev,
  Gallier, Kaminer, and Kumar}]{Keller2014}
Keller, J., Thakur, D., Dobrokhodov, V., Jones, K., Likhachev, M., Gallier, J.,
  Kaminer, I., and Kumar, V., \enquote{Coordinated commencement of pre-planned
  routes for fixed-wing UAS starting from arbitrary locations-a near real-time
  solution,} \emph{2014 International Conference on Unmanned Aircraft Systems,
  ICUAS 2014 - Conference Proceedings}, IEEE Computer Society, 2014, pp.
  552--561.
\newblock \doi{10.1109/ICUAS.2014.6842297}.

\bibitem[{Kallmann(2010)}]{Kallmann2010a}
Kallmann, M., \enquote{Shortest Paths with Arbitrary Clearance from Navigation
  Meshes,} \emph{Proceedings of the 2010 ACM SIGGRAPH/Eurographics Symposium on
  Computer Animation}, Eurographics Association, Goslar, DEU, 2010, p.
  159–168.
\newblock \doi{10.2312/SCA/SCA10/159-168}.

\bibitem[{Zollars et~al.(2018{\natexlab{a}})Zollars, Cobb, and
  Grymin}]{Zollars2017b}
Zollars, M.~D., Cobb, R.~G., and Grymin, D.~J., \enquote{Simplex Optimal
  Control Methods for Urban Environment Path Planning,} \emph{AIAA Sci-Tech
  Information Systems Conference}, 2018{\natexlab{a}}, p.~16.
\newblock \doi{10.2514/6.2018-2259}.

\bibitem[{Zollars et~al.(2017)Zollars, Cobb, and Grymin}]{Zollars2017}
Zollars, M.~D., Cobb, R.~G., and Grymin, D.~J., \enquote{Simplex Solutions for
  Optimal Control Flight Paths in Urban Environments,} \emph{Journal of
  Aeronautics and Aerospace Engineering}, Vol.~6, 2017, p.~8.
\newblock \doi{10.4172/2168-9792.1000197}.

\bibitem[{Smith(2014)}]{Smith2014}
Smith, N.~E., \enquote{Optimal Collision Avoidance Trajectories for
  Unmanned/Remotely Piloted Aircraft,} Ph.D. thesis, Air Force Institute of
  Technology, 2014.

\bibitem[{Suplisson(2015)}]{Suplisson2015}
Suplisson, A.~W., \enquote{Optimal Recovery Trajectories for Automatic Ground
  Collision Avoidance Systems (Auto GCAS),} Ph.D. thesis, Air Force Institute
  of Technology, 2015.

\bibitem[{Betts(2008)}]{Betts2010}
Betts, J.~T., \enquote{Practical methods for optimal control and estimation
  using nonlinear programming,} \emph{Advances in design and control}, 2008, p.
  434.

\bibitem[{Chew(1987)}]{Chew1987}
Chew, L.~P., \enquote{Constrained Delaunay Triangulations,} \emph{Proceedings
  of the third annual symposium on computational geometry}, 1987, pp. 215--222.
\newblock \doi{10.1145/41958.41981}.

\bibitem[{Pearl~Dechter(1985)}]{Dechter1985}
Pearl~Dechter, R.~J., \enquote{Generalized best-first search strategies and the
  optimality of A*,} \emph{Journal of the ACM}, Vol.~32, 1985, pp. 505--536.
\newblock \doi{10.1145/3828.3830}.

\bibitem[{Lee and Preparata(1984)}]{Lee1984}
Lee, D.~T., and Preparata, F.~P., \enquote{Euclidean shortest paths in the
  presence of rectilinear barriers,} \emph{Networks}, Vol.~14, No.~3, 1984, pp.
  393--410.
\newblock \doi{https://doi.org/10.1002/net.3230140304},
  \urlprefix\url{https://onlinelibrary.wiley.com/doi/abs/10.1002/net.3230140304}.

\bibitem[{Chazelle(1982)}]{Chazelle1982}
Chazelle, B., \enquote{A theorem on polygon cutting with applications,}
  \emph{23rd IEEE Symposium on Foundations of Computer Science}, 1982, pp.
  339--349.
\newblock \doi{10.1109/SFCS.1982.58}.

\bibitem[{Hershberger and Snoeyink(1993)}]{Hershberger1994}
Hershberger, J., and Snoeyink, J., \enquote{Computing minimum length paths of a
  given homotopy class,} \emph{Computational Geometry: Theory and
  Applications}, Vol.~4, 1993, pp. 63--97.
\newblock \doi{10.1016/0925-7721(94)90010-8}.

\bibitem[{Zollars et~al.(2018{\natexlab{b}})Zollars, Cobb, and
  Grymin}]{Zollars2017c}
Zollars, M.~D., Cobb, R.~G., and Grymin, D.~J., \enquote{Optimal Path Planning
  for SUAS Waypoint Following in Urban Enviornments,} \emph{IEEE Aerospace
  Conference}, 2018{\natexlab{b}}, p.~10.
\newblock \doi{10.1109/AERO.2018.8396483}.

\end{thebibliography}
\end{document}